\documentclass[11pt]{article}
\usepackage{psfig}
\usepackage{amsmath}
\usepackage{geometry}                
\usepackage{graphicx}
\usepackage{amssymb}
\usepackage{epstopdf}

\usepackage {natbib}
\usepackage{epsfig}
\usepackage{url}

\begin{document}

\fontsize{11}{14.5pt}\selectfont

\begin{center}

{\small Technical Report No.\ 0707,
 Department of Statistics, University of Toronto}

\vspace*{0.9in}

{\Large\bf Nonlinear Models Using Dirichlet Process 
           Mixtures}\\[16pt]

 \parbox[t]{7.5cm} {
 \begin{center}
{\large Babak Shahbaba}\\[2pt]
 Dept. of Public Health Sciences, Biostatistics\\ 
 University of Toronto, \\
 Toronto, Ontario, Canada \\
 \texttt{babak@stat.utoronto.ca}\\[10pt]
 \end{center}
}
\hfill
\parbox[t]{8.6cm} {
\begin{center}
{\large Radford M. Neal}\\[2pt]
 Dept. of Statistics and Dept. of Computer Science,\\
 University of Toronto, \\
 Toronto, Ontario, Canada \\
 \texttt{radford@stat.utoronto.ca}\\[10pt]
 \end{center}
 }

9 March 2007

\end{center}

\vspace*{8pt}

\noindent \textbf{Abstract.} We introduce a new nonlinear model for classification, in which we model the joint distribution of response variable, $y$, and covariates, $x$, non-parametrically using Dirichlet process mixtures. We keep the relationship between $y$ and $x$ linear within each component of the mixture. The overall relationship becomes nonlinear if the mixture contains more than one component. We use simulated data to compare the performance of this new approach to a simple multinomial logit (MNL) model, an MNL model with quadratic terms, and a decision tree model. We also evaluate our approach on a protein fold classification problem, and find that our model provides substantial improvement over previous methods, which were based on Neural Networks (NN) and Support Vector Machines (SVM). Folding classes of protein have a hierarchical structure. We extend our method to classification problems where a class hierarchy is available. We find that using the prior information regarding the hierarchical structure of protein folds can result in higher predictive accuracy.  

\section{\hspace*{-7pt}Introduction}\vspace*{-12pt}
In regression and classification models, estimation of parameters and interpretation of results are easier if we assume a simple distributional form (e.g., normality) and regard the relationship between response variable and covariates as linear. However, the performance of the model obtained depends on the appropriateness of these assumptions. Poor performance may result from assuming wrong distributions, or regarding relationships as linear when they are not. In this paper, we introduce a new model based on a Dirichlet process mixture of simple distributions, which is more flexible to capture nonlinear relationships. 

A Dirichlet process, $\mathcal{D}(G_0, \gamma)$, with baseline distribution $G_0$ and scale parameter $\gamma$, is a distribution over distributions. \cite{ferguson73} introduced the Dirichlet process as a class of prior distributions for which the support is large, and the posterior distribution is manageable analytically. Using the Polya urn scheme, \cite{blackwell73} showed that the distributions sampled from a Dirichlet process are discrete almost surely. The idea of using a Dirichlet process as the prior for the mixing proportions of a simple distribution (e.g., Gaussian) was first introduced by \cite{antoniak74}. 

We will describe the Dirichlet process mixture model as a limit of finite mixture model (see \cite{neal00} for further description). Suppose $y_1, ..., y_n$ are drawn independently from some unknown distribution. We can model the distribution of $y$ as a mixture of simple distributions such that:
\begin{eqnarray*}
P(y) & = & \sum_{c=1}^{C}p_c f(y|\phi_c)
\end{eqnarray*}
Here, $p_c$ are the mixing proportions, and $f$ is a simple class of distributions, such as normal with $\phi = (\mu, \sigma)$. We first assume that the number of mixing components, $C$, is finite. In this case, a common prior for $p_c$ is a symmetric Dirichlet distribution:
\begin{eqnarray*}
P(p_1, ..., p_C) & = & \frac{\Gamma(\gamma)}{\Gamma{(\gamma/C)}^{C}}\prod_{c=1}^{C}p_{c}^{(\gamma/C )-1}
\end{eqnarray*}
where $p_c \ge 0$ and $\sum p_c = 1$. Parameters $\phi_c$ are assumed to be independent under the prior with distribution $G_0$. We can use mixture identifiers, $c_i$, and represent the above mixture model as follows  \citep{neal00}: 

\parbox{14cm}{
\begin{eqnarray*}\label{eq:finiteMix}
y_i | c_i, \phi & \sim & F(\phi_{c_i}) \\
c_i | p_1, ..., p_C & \sim & Discrete(p_1, ..., p_C) \\
p_1, ..., p_C & \sim & Dirichlet(\gamma/C, ...., \gamma/C) \\
\phi_{c}& \sim & G_0 
\end{eqnarray*}
}
\hfil
\parbox{1cm}{\begin{equation}\end{equation}}
By integrating over the Dirichlet prior, we can eliminate mixing proportions, $p_c$, and obtain the following conditional distribution for $c_i$:
\begin{eqnarray} \label{eq:condProbFinite}
P(c_i = c | c_1, ..., c_{i-1}) & = & \frac{n_{ic}+\gamma/C }{i-1+\gamma}
\end{eqnarray}
Here, $n_{ic}$ represents the number of data points previously (i.e., before the $i^{th}$) assigned to component $c$. The probability of assigning each component to the first data point is $1/C$. As we proceed, this probability becomes higher for components with larger numbers of samples (i.e., larger $n_{ic}$). 

When $C$ goes to infinity, the conditional probabilities (\ref{eq:condProbFinite}) reach the following limits: 

\parbox{15cm}{
\begin{eqnarray*}
P(c_i = c | c_1, ..., c_{i-1}) & \to & \frac{n_{ic} }{i-1+\gamma}\\
P(c_i \ne c_j \forall j < i | c_1, ..., c_{i-1}) & \to & \frac{\gamma}{i-1+\gamma}
\end{eqnarray*}
}
\hfil
\parbox{1cm}{\begin{equation}\end{equation}}   
As a result, the conditional probability for $\theta_i$, where $\theta_i = \phi_{c_i}$, becomes
\begin{eqnarray}
\theta_i | \theta_1, ..., \theta_{i-1} & \sim & \frac{1}{i-1+\gamma}\sum_{j<i}\delta(\theta_j) + \frac{\gamma}{i-1+\gamma}G_0
\end{eqnarray}
where $\delta(\theta)$ is a point mass distribution at $\theta$. This is equivalent to the conditional probabilities implied by the Dirichlet process mixture model, which has the following form:
\begin{eqnarray} \label{eq:dp}
y_i | \theta_i & \sim & F(\theta_i) \nonumber \\
\theta_i | G & \sim & G \\
G & \sim & \mathcal{D}(G_0, \gamma) \nonumber
\end{eqnarray}
That is, the limit of the finite mixture model (\ref{eq:finiteMix}) is equivalent to the Dirichlet process mixture model (\ref{eq:dp}) as the number of components goes to infinity. $G$ is the distribution over $\theta$'s, and has a Dirichlet process prior, $\mathcal{D}$. The parameters of the Dirichlet process prior are $G_0$, a distribution from which $\theta$'s are sampled, and $\gamma$, a positive scale parameter that controls the number of components in the mixture, such that a larger $\gamma$ results in a larger number of components. Phrased this way, each data point, $i$, has its own parameters, $\theta_i$, drawn from a distribution that is drawn from a Dirichlet process prior. But since distributions drawn from a Dirichlet process are discrete (almost surely), the $\theta_i$ for different data points may be the same.  

\cite{bush96}, \cite{escobar95}, \cite{maceachern98}, and \cite{neal00} have used this method for density estimation. \cite{muller96} used Dirichlet process mixtures for curve fitting. They model the joint distribution of data pairs $(x_i, y_i)$ as a Dirichlet process mixture of multivariate normals. The conditional distribution, $P(y|x)$, and the expected value, $E(y|x)$, are estimated based on this distribution for a grid of $x$'s (with interpolation) to obtain a nonparametric curve. The application of this approach \citep[as presented by][]{muller96} is restricted to continuous variables. Moreover, this model is feasible only for problems with a small number of covariates, $p$. For data with moderate to large dimensionality, estimation of the joint distribution is very difficult both statistically and computationally. This is mostly due to the difficulties that arise when simulating from the posterior distribution of large full covariance matrices. In this approach, if a mixture model has $C$ components, the set of full covariance matrices have $Cp(p+1)/2$ parameters. For large $p$, the computational burden of estimating these parameters might be overwhelming. Estimating full covariance matrices can also cause statistical difficulties since we need to assure that covariance matrices are positive semidefinite. Conjugate priors based the inverse Wishart distribution satisfy this requirement, but they lack flexibility \citep{daniels99}. Flat priors may not be suitable either, since they can lead to improper posterior distributions, and they can be unintentionally informative \citep{daniels99}. A common approach to address these issues is to use decomposition methods in specifying priors for full covariance matrices \citep[see for example,][]{daniels99, cai06}. Although this approach has demonstrated some computational advantages over direct estimation of full covariance matrices, it is not yet feasible for high-dimensional variables. For example, \cite{cai06} recommend their approach only for problems with less than 20 covariates. 

We introduce a new nonlinear Bayesian model, which also non-parametrically estimates the joint distribution of the response variable, $y$, and covariates, $x$, using Dirichlet process mixtures. Within each component, we assume the covariates are independent, and model the dependence between $y$ and $x$ using a linear model. Therefore, unlike the method of \cite{muller96}, our approach can be used for modeling data with a large number of covariates, since the covariance matrix for one mixture component is highly restricted. Moreover, this method can be used for categorical as well as continuous response variables by using a generalized linear model instead of the linear model of each component. 

Our focus in this paper is on classification models with a multi-category response. We also show how our method can be extended to classification problems where classes have a hierarchical structure, and to problems with multiple sources of information. The next section describes our methodology. In Section 3, we illustrate our approach and evaluate its performance based on simulated data. In Section 4, we present the results of applying our model to an actual classification problem, which attempts to identify the folding class of a protein sequence based on the composition of its amino acids. Folding classes of protein have a hierarchical structure. In Section 5, we extend our approach to classification problems of this sort where a class hierarchy is available, and evaluate the performance of this new model on the protein fold recognition dataset. Section 6 shows how this approach can be used for multiple sources of information. Finally, Section 7 is devoted to discussion, future directions and limitations of the proposed method. 

\section{\hspace*{-7pt} Methodology}\vspace*{-12pt}
Consider a classification problem with continuous covariates, $x=(x_1, ..., x_p)$, and a categorical response variable, $y$, with $J$ classes. To model the relationship between $y$ and $x$, we model the joint distribution of $y$ and $x$ non-parametrically using Dirichlet process mixtures. Within each component of the mixture, the relationship between $y$ and $x$ is assumed to be linear. The overall relationship becomes nonlinear if the mixture contains more than one component. This way, while we relax the assumption of linearity, the flexibility of the relationship is controlled. Our model has the following form:
\begin{eqnarray*}
y_i, x_{i1}, ..., x_{ip} | \theta_i & \sim & F(\theta_i) \\
\theta_i | G & \sim & G \\
G & \sim & \mathcal{D}(G_0, \gamma)
\end{eqnarray*}
where $i=1, ..., n$ indexes the observations, and $l = 1, ..., p$ indexes the covariates.  In our model, $\theta = (\mu, \sigma, \alpha, \boldsymbol{\beta})$, and the component distributions, $F(\theta)$, are defined based on $P(y, x) = P(x)P(y|x)$ as follows:
\begin{eqnarray*}
x_{il} & \sim & N(\mu_{l}, \sigma_{l}^{2})\\
P(y_i=j|x_i, \boldsymbol{\alpha}, \boldsymbol{\beta}) & = &\frac{\exp(\alpha_{j} + x_i\boldsymbol{\beta}_{j})}{\sum_{j'=1}^{J} \exp(\alpha_{j'}+x_i\boldsymbol{ \beta}_{j'})}
\end{eqnarray*}
Here, the parameters $\mu=(\mu_1, ..., \mu_p)$ and $\sigma = (\sigma_1, ..., \sigma_p)$ are the means and standard deviations of covariates in each component. The component index, $c$, is omitted for simplicity. Within a component, $\alpha = (\alpha_1, ..., \alpha_J)$, and $\boldsymbol{\beta }= (\boldsymbol{\beta_1, ..., \beta_J})$ are the parameters of the multinomial logit (MNL) model, and $J$ is the number of classes. The entire set of regression coefficients, $\boldsymbol{\beta }$, can be presented as a $p \times J$ matrix. This representation is redundant, since one of the $\boldsymbol{\beta}_j$'s (where $j = 1, ..., J$) can be set to zero without changing the set of relationships expressible with the model, but removing this redundancy would make it difficult to specify a prior that treats all classes symmetrically. In this parameterization, what matters is the difference between the parameters of different classes. 

Although the covariates in each component are assumed to be independent with normal priors, this independence of covariates exists only locally (within a component). Their global (over all components) dependency is modeled by assigning data to different components (i.e., clustering). The relationship between $y$ and $x$ within a component is captured using an MNL model. Therefore, the relationship is linear locally, but nonlinear globally. 

We could assume that $y$ and $x$ are independent within components, and capture the dependence between the response and the covariates by clustering too. However, this may lead to poor performance (e.g., when predicting the response for new observations) if the dependence of $y$ on $x$ is difficult to capture using clustering alone. Alternatively, we could also assume that the covariates are dependent within a component. For continuous response variables, this becomes equivalent to the model proposed by \cite{muller96}. However, as we discussed above, this approach may be practically infeasible for problems with a moderate to large number of covariates. We believe that our method is an appropriate compromise between these two alternatives.

We define $G_0$ as follows:
\begin{eqnarray*}
\mu_{l}|\mu_{0}, \sigma_{0} & \sim & N(\mu_{0}, \sigma_{0}^{2}) \\
\log(\sigma^{2}_{l})|M_{\sigma}, V_{\sigma} & \sim & N(M_{\sigma}, V_{\sigma}^{2})\\
\alpha_j | \tau & \sim & N(0, \tau^{2}) \\
\beta_{jl} | \nu & \sim & N({0}, \nu^{2})
\end{eqnarray*}
The parameters of $G_0$ may in turn depend on higher level hyperparameters. For example, we can regard the variances of coefficients as hyperparameters with the following priors:
\begin{eqnarray*}
\log(\tau^{2}) | M_{\tau}, V_{\tau} & \sim & N(M_{\tau}, V_{\tau}^{2}) \\
\log(\nu^{2}) | M_{\nu}, V_{\nu} & \sim & N(M_{\nu}, V_{\nu}^{2})
\end{eqnarray*}

We use MCMC algorithms for posterior sampling. Samples simulated from the posterior distribution are used to estimate posterior predictive probabilities. We predict the response values for new cases based on these probabilities. For a new case with covariates $x'$, the posterior predictive probability of response variable, $y'$, is estimated as follows:
\begin{eqnarray*}
P(y'=j|x') = \frac{P(y'=j, x')}{P(x')}
\end{eqnarray*}
where 
\begin{eqnarray*}
{P(y'=j, x')} & = & {\frac{1}{S}\sum_{s=1}^{S} P(y'=j, x'|G_0, \theta^{(s)})}\\
P(x') & = & {\frac{1}{S}\sum_{s=1}^{S} P(x'|G_0, \theta^{(s)})}
\end{eqnarray*}
Here, $S$ is the number of post-convergence samples from MCMC, and $\theta^{(s)}$ represents the set of parameters obtained at iteration $s$.

\cite{neal00} presented several possible algorithms for sampling from the posterior distribution of Dirichlet process mixtures. In this research, we use Gibbs sampling with auxiliary parameters (Neal's algorithm 8). This approach is similar to the algorithm proposed by \cite{maceachern98}, with a difference that the auxiliary parameters exist only temporarily. To improve the MCMC sampling, after each update using auxiliary variables, we update the component parameters using their corresponding data points. For a complete description of this method, see the paper by \cite{neal00}. All our models are coded in MATLAB and are available online at \url{http://www.utstat.utoronto.ca/~babak}.  

\begin{figure}[t] 
\begin{center} 
\epsfig{file=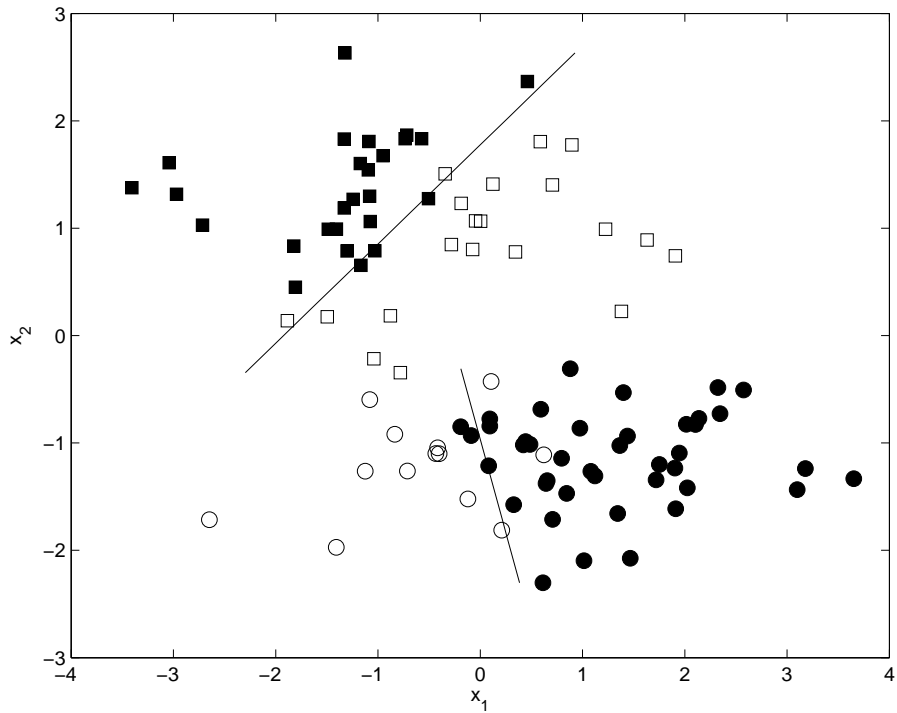, width=.85\textwidth}
\caption{An illustration of our model for a binary (black and white) classification problem with two covariates. Here, the mixture has two components, which are shown with circles and squares. In each component, an MNL model separates the two classes into ``black'' or ``white'' with a linear decision boundary.} %
 \label{dpExample1} %
\end{center} 
\end{figure}

In Figure \ref{dpExample1}, we show a state from an MCMC simulation for our model in which there are two covariates and the response variable is binary. In this iteration, our model has identified two components (circles and squares). Within a component, two classes (stars and crosses) are separated using an MNL model. Note, the decision boundaries shown are component specific. The overall decision boundary, which is a smooth function, is not shown in this figure. In our approach, division of the data into components and fitting of MNL models are performed simultaneously. 

\section{\hspace*{-7pt}Results for synthetic data}\vspace*{-12pt}
In this section, we illustrate our approach, henceforth called dpMNL, using synthetic data. We compare our model to a simple MNL model, an MNL model with quadratic terms (i.e., $x_{l}x_{k}$, where $l=1, ..., p$ and $k=1, ..., p$), referred to as qMNL, and a decision tree model \citep{breiman93} that uses 10-fold cross-validation for pruning. For the simple MNL model, we use both Bayesian and maximum likelihood estimation. The models are compared with respect to their accuracy rate and the $F_1$ measure. Accuracy rate is defined as the percentage of the times the correct class is predicted. $F_1$ is a common measurement in machine learning and is defined as:
\begin{eqnarray*}
F_1 &= & \frac{1}{J} \sum_{j=1}^{J}\frac{2A_j}{2A_j+B_j+C_j}
\end{eqnarray*}
where $A_j$ is the number of cases which are correctly assigned to class $j$, $B_j$ is the number cases incorrectly assigned to class $j$, and $C_j$ is the number of cases which belong to the class $j$ but are assigned to other classes. 

We do two tests. In the first test, we generate data according to the dpMNL model. Our objective is to evaluate the performance of our model when the distribution of data is comprised of multiple components. In the second test, we generate data using a smooth nonlinear function. Our goal is to evaluate the robustness of our model when data actually come from a different model.

For the first test, we compare the models using a synthetic four-way classification problem with 5 covariates. Data are generated according to our model with $G_0$ being the following prior:
\begin{eqnarray*}
\mu_{l} & \sim & N(0, 1) \\
\log(\sigma^{2}_{l}) & \sim & N(0 , 2^2)\\
\log(\tau^{2}) & \sim & N(0 , 0.1^{2}) \\
\log(\nu^{2}) & \sim & N(0, 2^2)
\end{eqnarray*}
Note that $\alpha_j | \tau \sim N(0, \tau^{2})$, and  $\beta_{jl} | \nu \sim N(0, \nu^{2})$, where $l = 1, ..., 5$ and $j = 1, ..., 4$. From the above baseline prior, we sample two components, $\theta_1$ and $\theta_2$, where $\theta = (\mu, \sigma, \eta, \nu, \alpha, \boldsymbol{\beta})$. For each $\theta$, we generate 5000 data points by first drawing $x_{il} \sim N(\mu_{l}, \sigma_{l})$ and then sampling $y$ using the following MNL model:
\begin{eqnarray*}
P(y=j|x, \boldsymbol{\alpha}, \boldsymbol{\beta}) & = &\frac{\exp(\alpha_{j} + x\boldsymbol{\beta}_{j})}{\sum_{j'=1}^{J} \exp(\alpha_{j'}+x\boldsymbol{ \beta}_{j'})}
\end{eqnarray*}
The overall sample size is $10000$. We randomly split the data to the training set, with 100 data points, and test set, with 9900 data points. We use the training set to fit the models, and use the independent test set to evaluate their performance. The regression parameters of the Bayesian MNL model with Bayesian estimation and the qMNL model have the following priors:
\begin{eqnarray*}
\alpha_j | \tau & \sim & N(0, \tau^{2}) \\
\beta_{jl} | \nu & \sim & N({0}, \nu^{2})  \\
\log(\eta) & \sim & N(0, 1^2) \\
\log(\nu) & \sim & N(0, 2^2)
\end{eqnarray*} 

To fit the decision tree models \citep{breiman93}, we used the available functions in MATLAB. These functions are treefit, treetest (for cross-validation) and treeprune.  

The above procedure was repeated 50 times. Each time, new $\theta_1$ and $\theta_2$ were sampled from the prior, and a new dataset was created based on these $\theta$'s. We used Hamiltonian dynamics \citep{neal93} for updating the regression parameters, $\alpha$'s and $ \boldsymbol{\beta}$'s. For all other parameters, we used single-variable slice sampling \citep{neal03} with the ``stepping out'' procedure to find an interval around the current point, and then the ``shrinkage'' procedure to sample from this interval. We also used slice sampling for updating the concentration parameter $\gamma$, where $\log(\gamma)\sim N(-3, 2^2)$. This prior encourages smaller values of $\gamma$, which results in smaller number of components. Note that the likelihood for $\gamma$ depends only on $C$, the number of unique components \citep{neal00, escobar95}. For all models we ran 5000 MCMC iterations to sample from the posterior distributions. We discarded the initial 500 samples and used the rest for prediction.

The average results (over 50 repetitions) are presented in Table \ref{simResults1}. As we can see, our dpMNL model provides better results compared to all other models. The improvements are statistically significant ($p$-values $< 0.001$ based accuracy rates) using a paired \emph{t}-test with $n=50$.
\begin{table}
\begin{center}
\begin{tabular}{l || c |c|c}
Model &   Accuracy (\%) & $F_1$ (\%)\\
\hline\hline
Baseline & 45.57 & 15.48 \\
\hline
MNL (Maximum Likelihood) &  77.30 &  66.65  \\
\hline
MNL & 78.39 & 66.52 \\
\hline 
qMNL & 83.60 & 74.16 \\
\hline
Tree (Cross Validation) & 70.87& 55.82 \\
\hline
dpMNL  & \bfseries 89.21 & \bfseries 81.00
\end{tabular}
\end{center}
\caption{Simulation 1: the average performance of models based on 50 simulated datasets. The Baseline model assigns test cases to the class with the highest frequency in the training set.}
\label{simResults1}
\end{table}%

\begin{figure}[t] 
\begin{center} 
\epsfig{file=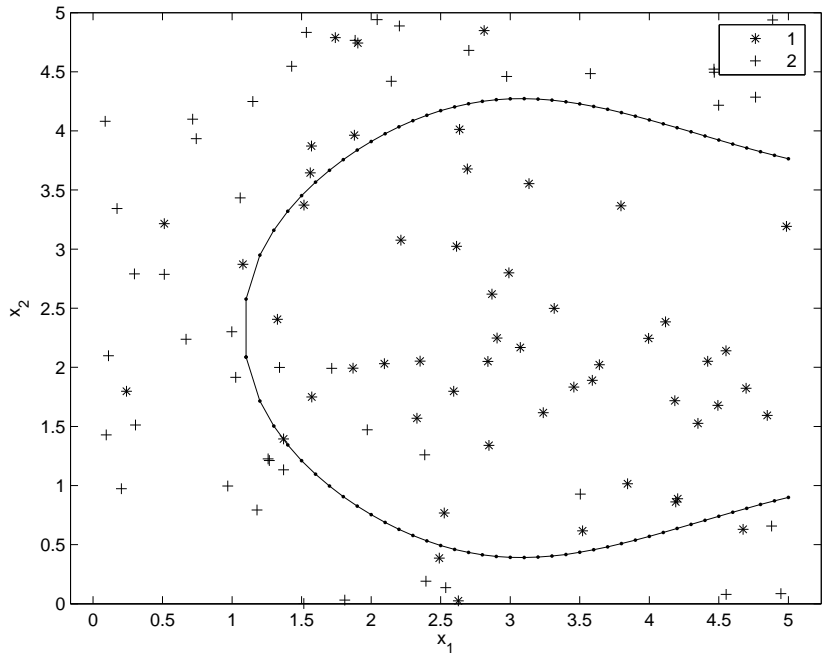, width=.78\textwidth}
\caption{A random sample generated according to Simulation 2 with $a_3=0$. The dotted line is the optimal boundary function.} %
 \label{dpExample2} %
\end{center} 
\end{figure}

Since the data were generated according to the dpMNL model, it is not surprising that this model had the best performance compared to other models. In fact, as we increase the number of components, the amount of improvement using our model becomes more and more substantial (results not shown). To evaluate the robustness of the dpMNL model, we performed another test. This time, we generated $x_{i1}, x_{i2}, x_{i3}$ (where $i=1, ..., 10000$) from the $Uniform(0, 5)$ distribution, and generated a binary response variable, $y_i$, according the following model:
\begin{eqnarray*}
 P(y=1| x) &= &\frac{1}{1+ \exp[a_1 \sin(x_{1}^{1.04}+1.2) + x_{1} \cos(a_2 x_{2} + 0.7) + a_3 x_{3} - 2]}
\end{eqnarray*}
where $a_1$, $a_2$ and $a_3$ are randomly sampled from $N(1, 0.5^2)$. The function used to generate $y$ is a smooth nonlinear function of covariates. The covariates are not clustered, so the generated data do not conform with the assumptions of our model. Moreover, this function includes a completely arbitrary set of constants to ensure the results are generalizable. Figure \ref{dpExample2} shows a random sample from this model for $a_3 = 0$. In this figure, the dotted line is the optimal decision boundary. 

We generated 50 datasets ($n=10000$) using the above model. Each time, we sampled new covariates, $x$, new constant values, $a_1, a_2, a_3$, and new response variable, $y$. As before, models were trained on 100 data points, and tested on the remaining samples. The average results over 50 datasets are presented in Table \ref{simResults2}. As before, the dpMNL model provides significantly (all $p$-values are smaller than 0.001) better performance compared to all other models. This time, however, the performance of the qMNL model is closer to the results from the dpMNL model. 
\begin{table}
\begin{center}
\begin{tabular}{l || c |c|c}
Model &   Accuracy (\%) & $F_1$ (\%)\\
\hline\hline
Baseline & 61.96 & 37.99 \\
\hline
MNL (Maximum Likelihood) & 73.58 &  68.33 \\
\hline
MNL & 73.58 & 67.92 \\
\hline 
qMNL & 75.60 & 70.12 \\
\hline
Tree (Cross Validation) & 73.47 & 66.94 \\
\hline
dpMNL  & \bfseries 77.80 & \bfseries 73.13\\
\end{tabular}
\end{center}
\caption{Simulation 2: the average performance of models based on 50 simulated datasets. The Baseline model assigns test cases to the class with the highest frequency in the training set.}
\label{simResults2}
\end{table}%

\section{\hspace*{-7pt}Results for protein fold classification}\vspace*{-12pt}
In this section, we consider the problem of predicting a protein's 3D structure (i.e., folding class) based on its sequence. For this problem, it is common to presume that the number of possible folds is fixed, and use a classification model to assign a protein to one of the folding classes. There are more than 600 folding patterns identified in the SCOP (Structural Classification of Proteins) database \citep{loconte00}. In this database, proteins are considered to have the same folding class if they have the same major secondary structure in the same arrangement with the same topological connections.

We apply our model to a protein fold recognition dataset provided by \cite{ding01}. The proteins in this dataset are obtained from the PDB\_select database \citep{hobohm92, hobohm94} such that two proteins have no more than 35\% of the sequence identity for aligned subsequences larger than 80 residues. Originally, the resulting dataset included 128 unique folds. However, \cite{ding01} selected only 27 most populated folds (311 proteins) for their analysis. They evaluated their models based on an independent sample (i.e., test set) obtained from PDB-40D \cite{loconte00}. PDB-40D contains the SCOP sequences with less than 40\% identity with each other. \cite{ding01} selected 383 representatives of the same 27 folds in the training set with no more than 35\% identity to the training sequences. The training and test datasets are available online at \url{http://crd.lbl.gov/~cding/protein/}. These datasets include the length of protein sequences, and 20 other covariates based on the percentage composition of different amino acids. For a detail description of data, see \cite{dubchak95}. 

\cite{ding01} trained several Support Vector Machines (SVM) with nonlinear kernel functions, and Neural Networks (NN) with different architecture on this dataset. They also tried different classification schemes, namely, one versus others (OvO), unique one versus others (uOvO), and all versus all (AvA). The details for these methods can be found in their paper. The performance of these models on the test set is presented in Table \ref{foldResults}. 

We first centered the covariates so they have mean 0. We trained our MNL and dpMNL on the training set, and evaluated their performance on the test set. For these models, we used similar priors as the ones used in the previous section. However, the hyperparameters for the variances of regression parameters are more elaborate. We used the following priors for the MNL model:
\begin{eqnarray*}
\alpha_j | \eta & \sim & N(0, \eta^2) \\
\log(\eta^{2}) & \sim & N(0, 2^2) \\
\beta_{jl} | \xi, \sigma_l & \sim & N({0}, \xi^2 \sigma_{l}^{2})  \\
\log(\xi^{2}) & \sim & N(0, 1) \\
\log(\sigma_{l}^{2}) & \sim & N(-3, 4^2) 
\end{eqnarray*} 
Here, one hyperparameter, $\sigma_l$, is used to control the variance of all coefficients, $\beta_{jl}$ (where $j=1, ..., J$), for covariate $x_l$. If a covariate is irrelevant, its hyperparameter will tend to be small, forcing the coefficients for that covariate to be near zero. This method is called Automatic Relevance Determination (ARD), and was suggested by \cite{neal96}. We also used another hyperparameter, $\xi$, to control the overall magnitude of all $\beta$'s. This way, $\sigma_l$ controls the relevance of covariate $x_l$ compared to other covariates, and $\xi$ controls the overall usefulness of all covariates in separating all classes. The standard deviation of $\beta_{jl}$ is therefore equal to $\xi \sigma_{l}$. 

We used the same scheme for the MNL models in dpMNL. Note that, in this model one $\sigma_l$ controls all $\beta_{jlc}$, where $j=1, ..., J$ indexes classes, and $c=1, ..., C$ indexes the unique components in the mixture. Therefore, the standard deviation of $\beta_{jlc}$ is  $\xi\sigma_{l}\nu_c $. Here, $\nu_c$ is specific to each component $c$, and controls the overall effect of coefficients in that component. That is, while $\sigma$ and $\xi$ are global hyperparameters common between all components, $\nu_c$ is a local hyperparameter within a component. Similarly, the standard deviation of intercepts, $\alpha_{jc}$ in component $c$ is $\eta \tau_c$. We used $N(0, 1)$ as the prior for $\nu_c$ and $\tau_c$. 

We also needed to specify priors for $\mu_l$ and $\sigma_l$, the mean and standard deviation of covariate $x_l$, where $l = 1, ..., p$. For these parameters, we used the following priors:  
\begin{eqnarray*}
\mu_{lc}|\mu_{0, l}, \sigma_{0, l} & \sim & N(\mu_{0, l}, \sigma_{0, l}^{2}) \\
\mu_{0, l} & \sim & N(0, 5^2) \\
\log(\sigma_{0, l}^{2}) & \sim & N(0, 2^2) \\
\log(\sigma^{2}_{lc})|M_{\sigma,l}, V_{\sigma,l} & \sim & N(M_{\sigma, l}, V_{\sigma, l}^{2})\\
M_{\sigma, l} & \sim & N(0, 1^2) \\
\log(V_{\sigma, l}^{2}) & \sim & N(0, 2^2)
\end{eqnarray*}
As we can see, the priors depend on higher level hyperparameters. This provides a more flexible scheme. If, for example, the components are not different with respect to covariate $x_l$, the corresponding variance, $\sigma^{2}_{0, l}$, becomes small, forcing $\mu_{lc}$ close to their overall mean, $\mu_{0, l}$.    

For each of our Bayesian models discussed in this section (and also in the following sections), we performed four simultaneous MCMC simulations each of size 10000. The chains have different starting values. We discarded the first 1000 samples from each chain and used the remaining samples for predictions. For this problem, running multiple chains results in faster and more efficient sampling. Simulating the Markov chain for 10 iterations took about half a minute for MNL, and about 3 minutes for dpMNL, using a MATLAB implementation on an UltraSPARC III machine. 

The results for MNL and dpMNL models are presented in Table \ref{foldResults}. As a benchmark, we also present the results for the SVM and NN models developed by \cite{ding01} on the exact same dataset. As we can see, our linear MNL model provides better accuracy rate compared to the SVM and NN models developed by \cite{ding01}. Our dpMNL model provides an additional improvement over the MNL model. This shows that there is in fact a nonlinear relationship between folding classes and the composition of amino acids, and our nonlinear model could successfully identify this relationship. 

It is worth noting the performance of the NN models is influenced by many design choices, and by model assumptions. We found that Bayesian neural networks model \citep{neal96} had better performance than the NN model of \cite{ding01}. Our NN model performs very similarly to the performance of the dpMNL model.

\begin{table}
\begin{center}
\begin{tabular}{l || c |c|c}
Model &   Accuracy (\%) & $F_1$ (\%)\\
\hline\hline
NN-OvO &  20.5 & -   \\
\hline
SVM-OvO &  43.5 & - \\
\hline
SVM-uOvO &  49.4 & - \\
\hline
SVM-AvA & 44.9 & - \\
\hline 
MNL & 50.0 & 41.2\\
\hline
dpMNL  & \bfseries 58.6 & \bfseries 53.0
\end{tabular}
\end{center}
\caption{Performance of models based on protein fold classification data. NN and SVM use maximum likelihood estimation, and are developed by \cite{ding01}.}
\label{foldResults}
\end{table}%

\section{\hspace*{-7pt}Extension to hierarchical classes}\vspace*{-12pt}
In the previous section, we modeled the folding classes as a set of unrelated entities. However, these classes are not completely unrelated, and can be grouped into four major structural classes known as $\alpha$, $\beta$, $\alpha / \beta$, and $\alpha+\beta$. \cite{ding01} show the corresponding hierarchical scheme (Table 1 in their paper). We have previously introduced a new approach for modeling hierarchical classes \citep{shahbabaBMC06, shahbabaBA07}. In this approach, we use a Bayesian form of the multinomial logit model, with a prior that introduces correlations between the parameters for classes that are nearby in the hierarchy. 

\begin{figure}
\unitlength .8mm
\begin{picture}(140,77)(-21,0)
\put(55,60){\makebox(0,0)[cc]{$\boldsymbol{\phi_{11}}$ }}

\put(110,60){\makebox(0,0)[cc]{$\boldsymbol{\phi_{12}}$ }}

\put(31.5,36){\makebox(0,0)[cc]{}}

\put(19,28){\makebox(0,0)[cc]{$\boldsymbol{\phi_{21}}$ }}

\put(42,28){\makebox(0,0)[cc]{$\boldsymbol{\phi_{22}}$ }}

\put(119,28){\makebox(0,0)[cc]{$\boldsymbol{\phi_{23}}$ }}

\put(142,28){\makebox(0,0)[cc]{$\boldsymbol{\phi_{24}}$ }}

\put(15,9){\makebox(0,0)[cc]{Class 1}}

\put(45,9){\makebox(0,0)[cc]{Class 2}}

\put(115,9){\makebox(0,0)[cc]{Class 3}}

\put(145,9){\makebox(0,0)[cc]{Class 4}}

\put(10,2){\makebox(0,0)[cc]{$\boldsymbol{\beta_1 = \phi_{11} + \phi_{21} }$}}

\put(51,2){\makebox(0,0)[cc]{$\boldsymbol{\beta_2 = \phi_{11} + \phi_{22} }$}}

\put(111,2){\makebox(0,0)[cc]{$\boldsymbol{\beta_3 = \phi_{12} + \phi_{23} }$}}

\put(152,2){\makebox(0,0)[cc]{$\boldsymbol{\beta_4 = \phi_{12} + \phi_{24} }$}}

\linethickness{0.3mm}
\multiput(15,15)(0.12,0.2){125}{\line(0,1){0.2}}
\linethickness{0.3mm}
\multiput(30,40)(0.12,-0.2){125}{\line(0,-1){0.2}}
\linethickness{0.3mm}
\linethickness{0.3mm}
\multiput(115,15)(0.12,0.2){125}{\line(0,1){0.2}}
\linethickness{0.3mm}
\multiput(130,40)(0.12,-0.2){125}{\line(0,-1){0.2}}
\linethickness{0.3mm}
\put(80,70){\circle{8.33}}

\put(80,70){\makebox(0,0)[cc]{1}}

\put(105,75){\makebox(0,0)[cc]{}}

\linethickness{0.3mm}
\multiput(30,50)(0.4,0.12){125}{\line(1,0){0.4}}
\linethickness{0.3mm}
\multiput(80,65)(0.4,-0.12){125}{\line(1,0){0.4}}
\linethickness{0.3mm}
\put(30,45){\circle{8.33}}

\linethickness{0.3mm}
\put(130,45){\circle{8.33}}

\put(30,45){\makebox(0,0)[cc]{2}}

\put(130,45){\makebox(0,0)[cc]{3}}

\end{picture}

\caption{A simple representation of our hierarchical classification model.}
\label{simpleHierarchy}
\end{figure}
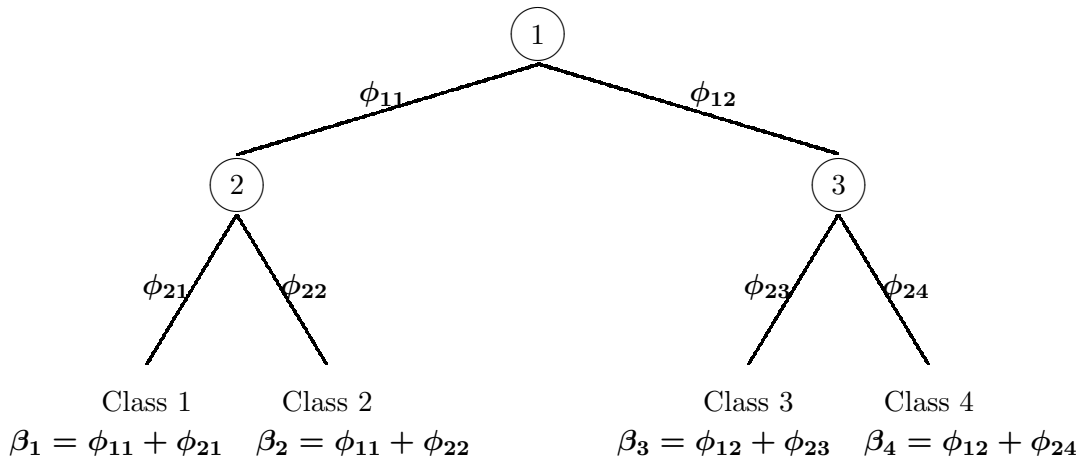

Figure \ref{simpleHierarchy} illustrates this approach using a simple hierarchical structure. For each branch in the hierarchy, we define a different set of parameters, $\phi$. Our model classifies objects to one of the end nodes using an MNL model whose regression coefficients for class $j$ are represented by the sum of the parameters for all the branches leading to that class. Sharing of common parameters (from common branches) introduces prior correlations between the parameters of nearby classes in the hierarchy. We refer to this model as corMNL.

In this section, we extend our nonlinear model to classification problems where classes have a hierarchical structure. For this purpose, we use a corMNL model, instead of MNL, to capture the relationship between the covariates, $x$, and the response variable, $y$, within each component. The results is a nonlinear model which takes the hierarchical structure of classes into account. We refer to this models as dpCorMNL.

Table \ref{foldResults2} presents the results for the two linear models (with and without hierarchy-base priors), and two nonlinear models (with and without hierarchy-based priors). In this table, ``parent accuracy'' refers to the accuracy of models based on the four major structural classes, namely $\alpha$, $\beta$, $\alpha / \beta$. When comparing the hierarchical models to their non-hierarchical counterparts, the advantage of using the hierarchy is apparent only for some measures (i.e., parent accuracy rate for corMNL, and the $F_1$ measure for dpCorMNL). As we can see, however, the dpCorMNL model provides a substantial improvement over corMNL. 

\begin{table}
\begin{center}
\begin{tabular}{l || c |c|c|c}
Model &   Accuracy (\%) & Parent accuracy (\%) & $F_1$ (\%)\\
\hline\hline
MNL & 50.0 & 76.5 & 41.2\\
\hline
corMNL & 49.5 & 77.9 & 41.4 \\
\hline
dpMNL  & 58.6 & \bfseries 79.9 & 53.0 \\
\hline
dpCorMNL  & \bfseries 59.1 & 79.4 & \bfseries 55.2
\end{tabular}
\end{center}
\caption{Comparison of hierarchical models (linear and nonlinear) with non-hierarchical models (linear and nonlinear) based on protein fold classification data.}
\label{foldResults2}
\end{table}%

\begin{table}[t]
\begin{center}
\begin{tabular}{l || c |c|c|c}
Model &   Accuracy (\%) & Parent accuracy (\%) & $F_1$ (\%)\\
\hline\hline
NN-OvO &  41.4 & - & - \\
\hline
SVM-OvO &  43.2 & - & - \\
\hline
SVM-uOvO &  49.4 & - & - \\
\hline
SVM-AvA & 56.5 & - & - \\
\hline
MNL & 56.5 & 80.4 & 51.4 \\
\hline
corMNL & 59.6 & 83.3 & 54.6\\
\hline
dpMNL & 60.4 & 82.0 & 55.9 \\
\hline
dpCorMNL  & \bfseries 61.4 & \bfseries 83.8 & \bfseries 57.8
\end{tabular}
\end{center}
\caption{Comparison of hierarchical models (linear and nonlinear) with non-hierarchical models (linear and nonlinear) based on protein fold classification data. The covariates are obtained from four different feature sets: composition of amino acids, predicted secondary structure, hydrophobicity, and normalized van der Waals volume.}
\label{foldResults3}
\end{table}%

\section{\hspace*{-7pt}Extension to multiple datasets}\vspace*{-12pt}
In order to improve the prediction of folding classes for proteins, \cite{ding01} combined the feature set based on amino acid compositions with 5 other feature sets, which were independently extracted based on various physico-chemical and structural properties of amino acids in the sequence. The additional features predicted secondary structure, hydrophobicity, normalized varn der Waals volume, polarity, and polarizability. Each data source has 21 covariates. For a detailed description of these features, see \cite{dubchak95}. \cite{ding01} added the above 5 datasets sequentially to the amino acid composition dataset. For prediction, they used a majority voting system, in which the votes obtained from models based on different features sets are combined, and the class with the most votes is regarded as predicted fold. Their results show that adding additional feature sets can improve the performance in some cases and can result in lower performance in some other cases. One main issue with this method is that it gives equal weights to votes based on different data sources. The underlying assumption, therefore, is that the quality of predictions is the same for all sources of information. This is, of course, not a realistic assumption for many real problems. In our previous paper \citep{shahbabaBMC06}, we provided a new scheme for combining different sources of information. In this approach, we use separate scale parameters, $\xi$, for each data source in order to adjust their relative weights automatically. This allows the coefficients from different sources of data to have appropriately different variances in the model.

For models developed by \cite{ding01}, the highest accuracy rate, $56.5$, was achieved only when they combined the covariates based on the composition of amino acids, secondary structure, hydrophobicity, and polarity. We also used these four datasets, and applied our models to the combined data. We used a different scale parameters, $\xi$, for each dataset. The results from our models are presented in Table \ref{foldResults3}. For comparison, we also present the results obtained by \cite{ding01} based on the same datasets. As we can see, this time, using the hierarchy results in more substantial improvements. Moreover, nonlinear models provided better performance compared to their corresponding linear models.

\section{\hspace*{-7pt}Conclusions and future directions}\vspace*{-12pt}

We introduced a new nonlinear classification model, which uses Dirichlet process mixtures to model the joint distribution of the response variable, $y$, and the covariates, $x$, non-parametrically. We compared our model to several linear and nonlinear alternative methods using both simulated and real data. We found that when the relationship between $y$ and $x$ is nonlinear, our approach provides substantial improvement over alternative methods. One advantage of this approach is that if the relationship is in fact linear, the model can easily reduce to a linear model by using only one component in the mixture. This way, it avoids overfitting, which is a common challenge in many nonlinear models.

We believe our model can provide more interpretable results. In many real problems, the identified components may correspond to a meaningful segmentation of data. Since the relationship between $y$ and $x$ remains linear in each segment, the results of our model can be expressed as a set of linear patterns for different segments of data. 

As mentioned above, for sampling from the posterior distribution, we used multiple chains which appeared to be sampling different regions of the posterior space. Ideally, we prefer to have one chain that can efficiently sample from the whole posterior distribution. In future, we intend to improve our MCMC sampling. For this purpose, we can use more efficient methods, such as the ``split-merge''  approach introduced by \cite{jain05} and the short-cut Metropolis method introduced by \cite{neal05}.

In this paper, we considered only continuous covariates. Our approach can be easily extended to situations where the covariate are categorical. For these problems, we need to replace the normal distribution in the baseline, $G_0$, with a more appropriate distribution. For example, when the covariate $x$ is binary, we can assume $x \sim Bernoulli(\mu)$, and specify an appropriate prior distribution (e.g., $Beta$ distribution) for $\mu$. Alternatively, we can use a continuous latent variable, $z$, such that $\mu = \exp(z)/ \{1+\exp(z)\}$. This way, we can still model the distribution of $z$ as a mixture of normals. For covariates with multinomial distribution, we can either extend the Bernoulli distribution by using ($\mu_1, ..., \mu_K$), where $K$ is the number of categories in $x$, or use $K$ continuous latent variables, $z_1, ..., z_K$, and set $\theta_j = \exp(z_j) / \sum_{j'}^{K}\exp(z_j')$. 

Our model can also be extended to problems where the response variable is not multinomial. For example, we can use this approach for regression problems with continuous response, $y$. The distribution of $y$ can be assumed normal within a component. We model the mean of this normal distribution as a linear function of covariates for cases that belong to that component. Other types of response variables (i.e., with Poisson distribution) can be handled in a similar way. 

Finally, our approach provides a convenient framework for semi-supervised learning, in which both labeled and unlabeled data are used in the learning process. In our approach, unlabeled data can contribute to modeling the distribution of covariates, $x$, while only labeled data are used to identify the dependence between $y$ and $x$. This is a quite useful approach for problems where the response variable is known for a limited number of cases, but a large amount of unlabeled data can be generated. One such problem is classification of web documents. In future, we will examine the application of our approach for these problems.


\end{document}